\documentclass[12, eqno]{amsart}
\usepackage{amsmath}
\usepackage{amscd,amsthm,amsfonts,amsopn,amssymb,verbatim}
\parskip =\medskipamount
\numberwithin{equation}{section}
\theoremstyle{remark}
\DeclareMathOperator{\supp}{supp\,}

\def\be{\begin{equation}}
\def\ee{\end{equation}}
\def\vp{\varphi}

\allowdisplaybreaks

\begin{document}

\title{Nonlinear Roth type theorems in finite fields}
\author{J.~Bourgain and M.-C.~Chang}
\address{J.~Bourgain, Institute for Advanced Study, Princeton, NJ 08540}
\email{bourgain@math.ias.edu}
\address{M.-C.~Chang, University of California, Riverside, CA 92521}
\email{mcc@math.ucr.edu}
\thanks{The authors were partially supported by NSF grants DMS-1301619 and DMS~1600154.}
\begin{abstract}
We obtain smoothing estimates for certain nonlinear convolution operators on prime fields, leading to quantitative nonlinear Roth type theorems.

Compared with the usual linear setting (i.e. arithmetic progressions), the nonlinear nature of the operators leads to different phenomena, both
qualitatively and quantitatively.
\end{abstract}
\maketitle

s
\section
{Introduction}
In this note, we study in the setting of prime fields certain nonlinear averages, motivated by the non-conventional ergodic averages considered in
particular in \cite {FK}.
While the results from \cite {FK} are of a more general nature, our emphasis are the quantitative aspects that also are in sharp contrast with the
linear case corresponding to arithmetic progressions.
As discussed below, this phenomenon leads to several questions related to generalization in various directions, worth further investigation.

For $f: \mathbb F_p\to \mathbb C$.
Denote
$$
\begin{aligned}
\mathbb E[f] &=\mathbb E_x[f]=\frac 1p \sum^{p-1}_{x=0} f(x)\\
\Vert f\Vert_r & =\Big(\frac 1p \sum_x |f(x)|^r\Big)^{\frac 1r}\\
\Vert f\Vert_{\ell^r} &= \Big(\sum_x | f(x)|^r\Big)^{\frac 1r}\\
\hat f(z)&= \frac 1p \sum_x e_p (-xz) f(x).
\end{aligned}
$$
Thus with this notation, Parseval's identity  reads thus
$$
\Vert f\Vert_2 = \Vert\hat f\Vert_{\ell^2}.
$$

Our first result is the following inequality for a certain nonlinear convolution.

\noindent
{\bf Theorem 1.1.}
{\sl Let $f_1, f_2:\mathbb F_p \to\mathbb C$. Define
\be\label{1.1}
F(x) =\frac 1p \sum_{y\in\mathbb F_p} f_1 (x+y) f_2 (x+y^2).\ee
Then
\be\label{1.2}
\Vert F-\mathbb E[f_1]\cdot \mathbb E[f_2]\Vert_2 \leq cp^{-\frac 1{10}} \Vert f_1\Vert_2 \cdot \Vert f_2\Vert_2
.\ee}

\noindent
{\bf Corollary 1.2.}
{\sl (quadratic Roth theorem on $\mathbb F_p$).

Let $A\subset \mathbb F_p$, $|A|=\delta p$ with $\delta >c_1 p^{-\frac 1{15}}$ and $c_1>0$ an appropriate constant.

Then there are $x\in A, y\in \mathbb F_p^*$ such that $x, x+y, x+y^2 \in A$

(and in fact $\gtrsim \delta^3p^2$ such triplets).}
\medskip

\noindent
{\bf Proof of Corollary 1.2.}

Let $f$ be the indicator function of $A$. Then $\Vert f\Vert_2=\delta^{\frac 12}$ and $\mathbb E(f)=\delta$.

Let $F$ be as in \eqref{1.1} with $f_1=f_2=f$. Then Cauchy-Schwarz and \eqref{1.2} imply
\be\label{1.3}
\begin{aligned}
\mathbb E_x \mathbb E_y [f(x)f(x+y) f(x+y^2)] =\;&\mathbb E_x[f F]\\
=\;&\mathbb E_x[f \cdot \mathbb E[f]^2]+\mathbb E_x[f \cdot (F-\mathbb E[f]^2)]\\
\geq
&\mathbb E[f]^3- \Vert f\Vert_2 \;\Vert F-\mathbb E[f]^2\Vert_2\\
\geq
&\mathbb E[f]^3-cp^{-\frac 1{10}} \Vert f\Vert_2^3\\\geq
&\delta^{3}-cp^{-\frac 1{10}}\delta^{\frac 32} \sim \delta^3.\quad\square
\end{aligned}
\ee

\noindent
{\bf Remarks}
\begin{itemize}
\item [(i)]
Recall Berend's construction [Be] providing a subset $A\subset \{ 0, 1, \ldots, p-1\}$ or $A\subset \mathbb
F_p$, $|A|>\delta p$ with $\delta\asymp e^{-\sqrt{\log p}}$ and containing no non-trivial triples $x, x+y,
x+2y\in A$.
Thus Corollary 1.2 in the  non-linear setting allows $\delta$ to be much smaller.

\item[(ii)] Theorem 1.1 is indeed a non-linear phenomenon, both qualitatively and quantitatively.
It was proven by N.~Frantzikinakis and B.~Kra \cite {FK} that if $(X, \mu, T)$ is totally ergodic measure
preserving probability system and $\{p_1(n), \ldots, p_k(n)\}$ an independent family of integer
polynomials, then for $f_1, \ldots, f_k\in L^\infty (\mu)$
\be\label{1.4}
\lim_{N\to\infty} \Big\Vert \frac 1N \sum^{N-1}_{n=0} f_1(T^{p_1(n)}x)\cdots f_k(T^{p_k(n)} x)-\prod^k_{i=1} \int f_i d\mu\Big\Vert_{L^2(\mu)}=0.\ee
The case $k=2, p_1(n)=n$, $p_2(n)=n^2$ (corresponding to Theorem 1.1) is due to Furstenberg and Weiss.
Thus what makes \eqref{1.4} a non-linear phenomenon (compared with the linear case $p_1(n) =n, p_2(n) =2n, \cdots,
p_k(n) =kn$, i.e. Szemeredi's theorem) is that the characteristic factor in the situation (1.4) turns out
to be trivial (this was first conjectured by V.~Bergelson).
On the quantitative side, the approach in \cite {FK} (based on the work of Host-Kra \cite{HK})
 does not seem to provide any estimates however, which is the main interest of Theorem 1.1.

\item[(iii)] The proof of Theorem 1.1 relies on Fourier analysis, much in the spirit of \cite{B}.
Compared with \cite {B}, certain simplification occurs due to the fact that oscillatory integrals are
replaced by certain exponential sums that allow simple (and optimal) bounds (by use of Weil's estimate).
A natural question is of course to what extent Theorem 1.1 generalizes to other non-linear settings
(see discussion at the end).
\end{itemize}

\section
{Proof of Theorem 1.1}

The argument follows closely the proof of Lemma 5 of \cite{B}.

Expanding $f_1, f_2$ in Fourier sum gives
\be\label{2.1}
F(x) =\mathbb E_y[f_1(x+y) f_2(x+y^2)] =\sum^{p-1}_{n_1, n_2=0} \hat f_1 (n_1) \hat f_2(n_2) c_{n_1, n_2
}
e_p\big((n_1+n_2)x\big)\ee
with (using quadratic Gauss sums evaluation)
\be\label{2.2}\begin{aligned}
c_{n_1, n_2} =&\;\;\;\mathbb E_y[e_p(n_1 y+n_2y^2)]\\=&\;\;\begin{cases}
\;\;\;\;\;\;\;\;\;\;\;\;1\;\;\;\;\;\;\;\;\;\;\;\;\;\;\;\;\;\;\;\;\;\;\;\text {    if } n_1=n_2=0\\ \;\;\;\;\;\;\;\;\;\;\;\;0\;\;\;\;\;\;\;\;\;\;\;\;\;\;\;\;\;\;\;\;\;\;\;\text {    if } n_2=0, n_1\not= 0\\
\frac1{\sqrt p} \big(\frac{n_2}{p}\big) e_p (-n_1^2 \, \overline {4n}_2)\sigma_p \;\;\text { if } n_2\not=0
\end{cases}\end{aligned}
\ee
and where
$$\sigma_p = \begin{cases} 1 \quad\text{ if } \;p\equiv 1\pmod 4\\ i \quad\text{ if } \;p\equiv 3\pmod 4\end{cases}\qquad\qquad\qquad$$
and $\bar x$ stands for the multiplicative inverse of $x\in\mathbb F_p^*$.

It follows that
$$
F=\mathbb E[f_1]\cdot \mathbb E[f_2] +\frac 1{\sqrt p}\,  \sum_{s=0}^{p-1} e_p(sx)\,  \sum_{n=1}^{p-1} \hat f_1(s-n)
\hat f_2 (n) \Big(\frac np\Big) K(s-n, n)\sigma_p,\qquad
$$and by Parseval,
\be\label{2.3}
\Vert F-\mathbb E[f_1]\cdot\mathbb E[f_2]\Vert_2 = \frac 1{\sqrt p} \Big\{\sum_s\Big|\sum_{n\not= 0}\hat f_1 (s-n)
\hat f_2 (n) \Big(\frac np\Big) K(s-n, n)\Big|^2 \Big\}^{\frac 12}\ee
where we denote
\be\label{2.4}
K(x, y) =\begin{cases}  e_p (-x^2 \,\overline {4y}\,) &\text { if } y\not =0\\ 0&\text{ otherwise}.
\end{cases}.\ee
Since we are aiming for an estimate of \eqref{2.3} in terms of $\Vert f_1\Vert_2 . \Vert f_2\Vert_2$, the
factor $\Big(\frac np\Big)$ may be absorbed in $\hat f_2(n)$.
Thus we need to analyze further
\be\label{2.5}
\Big(\sum_s\Big|\sum_{n\not= 0} \hat f_1(s-n) \hat f_2 (n) K(s-n, n)\Big|^2\Big)^{\frac 12}\ee
with $K$ given by \eqref{2.4}.

In order to bound \eqref{2.5} non-trivially (note that the trivial bound  by \hfill\break
$\sqrt p\;\Vert f_1\Vert_2\cdot \Vert
f_2\Vert_2$
is just insufficient and our aim is to gain an extra $p^{-\gamma}$ from additional cancellation exploiting
\eqref{2.4}), the first step is to invoke the following very general inequality, identical to Lemma 7 in \cite {B}
(we also use the notation $\int dx$ for $\sum_x$).

\noindent
{\bf Lemma 2.1.}
{\sl Let $f, g:\mathbb F_p \to\mathbb C$ and $K:\mathbb F_p\times\mathbb F_p \to\mathbb C$.
Then
\be\label{2.7}
\begin{aligned}
&\Big\Vert \int K(s-x, x) f(s-x) g(x) dx\Big\Vert_{\ell^2_s}=: I\qquad\\
\leqq \;\;&\Vert f\Vert_{\ell^2}^{1/2} \Vert f\Vert^{1/2}_{\ell^4} \Vert g\Vert_{\ell^2}\,\;
\cdot\\
&\;\;\;\;\;\;\Big \Vert \int K(x, s-x) \overline{K(x+u, s-x)} \, \overline{K(x, s+v-x)} \, K(x+u, s+v-x) dx\Big\Vert ^{1/4}_{\ell^2_{s, u, v}}
.
\end{aligned}
\ee}

We repeat the argument for selfcontainedness.

\begin{proof}
The proof follows from consecutive linearizations and applications of the Cauchy-Schwarz inequality.

Linearization of the $\ell_s^2$-norm yields
\be\label{2.8}
\begin{aligned}
I&\leqq \Vert g\Vert_2\,\Big\Vert\int K(s-x, x) f(s-x)\vp(s) ds\Big\Vert_{\ell^2_x} \quad
(\text {for some $\vp$ with } \Vert\vp\Vert_{\ell^2}=1)\\
&\leqq \Vert g\Vert_2 \,\Big\Vert \int K(s-x, x) \overline{K(s'-x, x)} f(s-x) \overline
{f(s'-x)}dx\Big\Vert^{1/2}_{\ell^2_{s, s'}}.\end{aligned}
\ee
Make the change of variable $x\leftrightarrow s-x $ and put $u=s'-s$. Thus
\be\label{2.9}
\begin{aligned}
&\;\Big\Vert\int K(s-x, x)\overline{K(s'-x, x)} f(s-x)\overline {f(s'-x)}dx\Big\Vert_{\ell_{s, s'}^2}\\
=&\;\Big\Vert \int K(x, s -x) \overline {K(x+u, s-x)} f(x)  \overline {f(x+u)}dx\Big\Vert_{\ell^2_{s, u}}\end{aligned}
\ee
Fixing the $u$-variable and again linearizing the $\ell^2_s$-norm gives
$$
\begin{aligned}
&\;\Big\Vert \int K(x, s-x)\overline {K(x+u, s-x)} f(x)  \, \overline {f(x+u)} dx\Big\Vert_{\ell_s^2}\\
\leqq &\;\Big\Vert\int K(x, s-x) \overline {K(x+u, s-x)} \;\overline {K(x, s+v-x)}\\
&\qquad\qquad\qquad\qquad\qquad K(x+u, s+v-x) dx\Big\Vert^{1/2}_{\ell^2_{s, v}}\Vert f\cdot \bar f_u\Vert_{\ell^2_x}.
\end{aligned}
$$
Hence, by Cauchy-Schwarz again expression \eqref{2.9} is bounded by
\be\label{2.10}
\begin{aligned}
\Big\Vert \int K(x, s -x)&\overline{K(x+u, s-x) } \, \overline {K(x, s+v-x)}\\
&\times K(x+u, s+v-x)dx\Big\Vert^{1/2}_{\ell^2_{s, u, v}} \Big\Vert\Vert f\cdot \bar
f_u\Vert_{\ell^2_x}\Big\Vert_{\ell_u^4}.\end{aligned}
\ee
Since the last factor in \eqref{2.10} is bounded by
$$
 \max_u \mathbb E_x [|f|^2 |f_u|^2]^{\frac 14}. \Vert f\Vert_{\ell^2}\leq \Vert f\Vert_{\ell^2}\Vert
f\Vert_{\ell^4}
$$
\eqref{2.7} follows from \eqref{2.8}-\eqref{2.10}.
This proves Lemma 2.1.
\end{proof}

Returning to \eqref{2.5}, let $M>0$ be a parameter and decompose
\be\label{2.11}\begin{aligned}
\hat f_1 (n)=&\;\;\hat f_{1, 0}(n)+\hat f_{1,1}(n), \;\;\text{ with } \\
\hat f_{1, 0}(n) =&\;\;\begin{cases} \;\hat f_1(n)\quad \text { if }\hat f_1(n) <\frac M{\sqrt p}
\Vert f_1\Vert_2\\\;0\quad\quad\quad \text{ otherwise. }\end{cases}.\end{aligned}\ee
Applying Lemma 2.6 with $f=\hat f_{1, 0}$ and $g=\hat f_2$ implies
\be\label{2.12}
\big(\sum_s\Big|\sum_{n\not= 0} \hat f_{1, 0} (s-n) \hat f_2(n) K(s-n, n)|^2\Big)^{\frac 12}\leq
\Vert \hat f_{1, 0} \Vert_{\ell^2}^{\frac 12} \;\Vert\hat f_{1, 0}\Vert_{\ell^4}^{\frac 12} \;\Vert \hat
f_2\Vert_{\ell^2}\,\Omega^{\frac 14},
\ee
where
\be\label{2.13}
\Omega=\Big\Vert \sum_x [K(x, s-x)\overline{K(x+u, s-x)} \, \overline {K(x, s+v-x)} K(x+u,
s+v-x)]\Big\Vert_{\ell^2_{s, u, v}}.\ee
Also $$ \Vert\hat f_{1, 0}\Vert_{\ell^2} \leq \Vert\hat f_1\Vert_{\ell^2} =\Vert f_1\Vert_2,\;\;\;
\Vert\hat f_2\Vert_{\ell^2} = \Vert f_2\Vert_2
$$ while by \eqref{2.11}
$$
\Vert\hat f_{1, 0}\Vert_{\ell_4} =\Big(\sum\Big|\hat f_{1, 0}(n)|^4\Big)^{\frac 14} \leq \Big(\frac M{\sqrt
p}\Big)^{\frac 12} \Vert f_1\Vert_2^{\frac 12}\Big(\sum|\hat f_1 (n)|^2\Big)^{\frac 14} =\Big(\frac M{\sqrt
p}\Big)^{\frac 12} \Vert f_1\Vert_2.
$$
Thus the right hand side of \eqref{2.12} is bounded by
\be\label{2.14}
\Big(\frac M{\sqrt p}\Big)^{\frac 14} \Vert f_1\Vert_2 \cdot \Vert f_2\Vert_2 \,\, \Omega^{\frac
14}.\ee
Contribution of $\hat f_{1, 1} $ in \eqref{2.5} is estimated trivially as follows
\be\label{2.15}
\Big(\sum_s\Big| \sum_{n\not= 0} \hat f_{1, 1}(s-n) \hat f_2 (n)K(s-n, n)\Big|^2\Big)^{\frac 12}\leq
\Big(\sum_s\Big(\sum_n |\hat f_{1, 1} (s-n)| \, |\hat f_2(n)|\Big)^2\Big)^{\frac 12}\ee
and noting that since $\hat f_{1, 1}(n)=0$ or $|\hat f_{1, 1}(n)|\geq \frac M{\sqrt p}\Vert f_1\Vert_2$, it follows
$\Vert f_1 \Vert^2_2 \big(\frac M{\sqrt p}\big)^2 |\supp \hat f_{1, 1}| \leq \Vert\hat f_{1,
1}\Vert^2_{\ell^2}\leq \Vert f_1\Vert_2^2$  implying $|\supp \hat f_{1, 1}|\leq \big(\frac {\sqrt p}M\big)^2$. By Cauchy-Schwarz
$$
\sum_n |\hat f_{1, 1} (s-n)| \, |\hat f_2(n)| \leq \frac {\sqrt p}M\Big(\sum_n|\hat f_{1, 1} (s-n)|^2 \,
|\hat f_2(n)|^2\Big)^{\frac 12}
$$
and \eqref{2.15} is bounded by
\be\label{2.16}
 \frac {\sqrt p}M \, \Vert \hat f_{1, 1}\Vert_{\ell^2} \cdot \Vert \hat f_2\Vert_{\ell^2}
\leq \frac {\sqrt p}M \, \Vert f_1\Vert_2\cdot \Vert f_2\Vert_2.
\ee

Adding \eqref{2.14}, \eqref{2.16} shows that
$$\begin{aligned}
 &\bigg(\sum_s\Big|\sum_{n\not= 0} \hat f_1(s-n) \hat f_2 (n) K(s-n, n)\Big|^2\bigg)^{\frac 12}
 \leq&\bigg[\bigg(\frac M{\sqrt p}\bigg)^{\frac 14} \, \Omega^{\frac 14} +\frac {\sqrt p}M\bigg]\,
\Vert f_1\Vert_2 \cdot \Vert f_2\Vert _2\end{aligned}
$$
and choosing $M$ appropriately, we have a bound on \eqref{2.5}
\be\label{2.17}
\bigg(\sum_s\Big|\sum_{n\not= 0} \hat f_1(s-n) \hat f_2 (n) K(s-n, n)\Big|^2\bigg)^{\frac 12} \leq \Omega^{\frac 15}\, \Vert f_1\Vert_2 \cdot \Vert f_2\Vert_2.\ee
Next, recalling (2.4), we will establish a non-trivial bound on (2.13), exploiting cancellation in $\mathbb
E_x$.

We have by (2.4), for $s\not= x, s+v\not= x$
$$
K(x, s-x) \overline{K(x+u, s-x)} = e_p \Big(-\frac {x^2}{4(s-x)}+\frac {(x+u)^2}{4(s-x)}\Big) =
e_p \Big( \frac {(2x+u) u}{4(s-x)}\Big).$$

Hence
\be\label{2.18}
\begin{aligned}
&\;\mathcal K_{x,s,u,v}\\:=&\;K(x, s-x)\overline{K(x+u, s-x)} \ \overline{ K(x, s+v-x)} \, K (x+u, s+v-x)\\ = &\;e_p \bigg(\frac{(2x+u)u}{4}
\Big(\frac 1{s-x}-\frac 1{s+v-x}\Big)\bigg)\\
=&\; e_p \bigg(\frac {uv(2x+u)}{(s-x)(s+v-x)}\bigg),\end{aligned}
\ee
and assuming $u, v\not= 0, u\not= -2s, -2s-2v$, by the following proposition, we have
\be\label{2.19}
\bigg|\sum_{\substack{ x=0\\x\not= s, s+v}}^{p-1} \mathcal K_{x,s,u,v}\bigg| \leq 3\sqrt p.\ee

The following estimate is from [Bom] (Theorem~5).

\noindent
{\bf Proposition 2.2.}
{\sl Let $f_1, f_2 \in \mathbb Z[X]$, $(f_1, f_2)=1$ and $\tilde f_1, \tilde f_2 \in\mathbb F_p[X]$ the
corresponding polynomials over $\mathbb F_p$, $\tilde f(x) =\frac {\tilde f_1(x)}{\tilde f_2(x)}$ where $x$
is to take only values with $p\nmid f_2(x)$. Define
\be\label{2.21}
S(\tilde f) =\sum_x e_p \big(\tilde f(x)\big). \ee
Then, assuming $\deg(\tilde f):= \deg (\tilde f_1) +\deg (\tilde f_2)\geq 1$, we have
\be\label{2.11}
|S(\tilde f)|\leq \big(n-2+\deg (\tilde f)_\infty\big) p^{\frac 12} +1\ee
with $n=$ the number of the poles and $(\tilde f)_\infty$ the divisor of the poles of $\tilde f$ over the
algebraic closure $\bar{\mathbb F}_p$ (including $\infty$ if necessary).}
\medskip

Applying Proposition 2.2 with $\tilde f(x)$ as in \eqref{2.18}, $n=2, \,\deg (\tilde f)_\infty=2$ so that indeed
$$
|S(f)| \leq 2 p^{\frac 12}+1< 3\sqrt p.
$$

If the assumptions in \eqref{2.19} are not fulfilled, estimate trivially $|\sum_x \mathcal K_{x,s,u,v}|\leq p$.
Using \eqref{2.19}, it follows that
\be\label{2.23}
\Omega \leq c\{p^3(\sqrt p)^2+p^2\cdot p^2\}^{\frac 12} \leq cp^2.\ee
Substituting in \eqref{2.17} gives
\be\label{2.24}
\Big(\sum_s\Big|\sum_{n\not= 0} \hat f_1(s-n) \hat f_2 (n) K(s-n, n)\Big|^2\Big)^{\frac 12} \leq cp^{2/5} \Vert f_1\Vert _2\, \Vert f_2\Vert_2.\ee
Hence by \eqref{2.3}
$$
\big\Vert F-\mathbb E[f_1]\cdot\mathbb E[f_2]\big\Vert_2  \leq cp^{-\frac 1{10}} \Vert f_1\Vert_2 \, \Vert f_2\Vert_2
$$which is \eqref{1.2}.
This proves Theorem 1.1.

\medskip

\noindent{\bf Remark.} One may wonder about the sharpness of inequality \eqref{1.2} in Theorem 1.1. Although we did not attempt to optimize our approach, it almost surely will not answer this presumably difficult question. We only note a few examples below.

\medskip

\noindent{\it Example 1.} Take $f_1(x)=f_2(x)=e_p(x)$. Then $F(x)=e_p(2x)\big(\frac 1p\sum_y e_p(y+y^2)\big)$, hence $|F(x)|=\frac 1{\sqrt p}$ and $\big\Vert F-\mathbb E[f_1]\;\mathbb E[f_2]\big\Vert_2=p^{-\frac 12}\Vert f_1\Vert_2\; \Vert f_2\Vert_2.$

\medskip

\noindent{\it Example 2.} For $i=1,2$, define
$$f_i(x)=\begin{cases}\;\;1\quad\quad\text{if } x=0\\\;\;0\quad\quad\text{otherwise}.\end{cases}$$
Thus $\;\mathbb E[f_i]=\frac 1p$, $\Vert f_i\Vert_2=p^{-\frac 12}$ and
$$F(x)=\begin{cases}\;\;\frac 1p\qquad\qquad \text{if }x=0, -1\\\;\;0\qquad\qquad\text{ otherwise}.\end{cases}$$
Therefore, $\big\Vert F-\mathbb E[f_1]\;\mathbb E[f_2]\big\Vert_2=O(p^{-\frac 32})
=O(p^{-\frac 12}\Vert f_1\Vert_2\; \Vert f_2\Vert_2)$ again.

\smallskip

The next example is slightly more interesting as it shows that \eqref{1.2} cannot hold with $p^{-\frac 1{10}}$ replaced by $p^{-\frac 1{2}}$.

\medskip

\noindent{\it Example 3.}  Recall the Fourier formulation \eqref{2.3} and \eqref{2.4}. The left hand side of \eqref{1.2} equals
$$\frac 1{\sqrt p} \Big\{\sum_s\Big|\sum_{n\not= 0}\hat f_1 (s-n)
\hat f_2 (n) \Big(\frac np\Big) e_p(-\overline{4}n)\,e_p(-\overline{4n}s^2)\Big|^2 \Big\}^{\frac 12}$$
and absorbing $\big(\frac np\big) e_p(-\overline{4}n)$ in the $\hat f_2 (n)$-factor, one obtains
\be\label{2.25}
\frac 1{\sqrt p} \Big\{\sum_s\Big|\sum_{n\not= 0}\hat f_1 (s-n)
\hat f_2 (n) \,e_p(-\overline{4n}s^2)\Big|^2 \Big\}^{\frac 12}.\ee

Next, we will define $\hat f_1$ and $\hat f_2$. Let $D$ be the product of the primes less than $\frac 1{10}\log p$ and $\mathcal D$ the set of divisors of $D$. Hence
$$D<p^{\frac 19} \;\; \text{ and } \;\; |\mathcal D| >\exp\bigg(c\;\frac{\log p}{\log\log p}\bigg).$$
Define
\be\label{2.26}
\hat f_1(x)=\begin{cases}\; (2D)^{-\frac 12} \qquad\quad \text{if }\; x\in 4D\cdot\{1,2,\ldots, 2D\}\\\;\;\;\; 0 \qquad\quad\quad \quad \text{otherwise}\end{cases}\ee
and
\be\label{2.27}
\hat f_2(x)=\begin{cases}\; |\mathcal D|^{-\frac 12} \qquad\quad \text{if }\; x\in 4D\cdot\mathcal D\\\;\;\;\; 0 \qquad\quad\quad\;\; \text{otherwise}.\end{cases}\qquad\qquad\;\;\;\;\;\ee
Hence $\Vert f_1\Vert_2=1=\Vert f_2\Vert_2$.

Setting $s=4Ds_1$ and $n=4Dn_1$, we have a lower bound on $\sqrt p\cdot$\eqref{2.25} as follows.
\be\label{2.28}\begin{aligned}
&\;\Big\{\sum_s\Big|\sum_{n\not= 0}\hat f_1 (s-n)
\hat f_2 (n) \,e_p(-\overline{4n}s^2)\Big|^2 \Big\}^{\frac 12}\\
\ge &\;\bigg\{\sum_{D<s_1\le 2D}\;\bigg|\sum_{n_1\in\mathcal D\cap(s_1-\{1,\ldots, 2D\})}(2D)^{-\frac 12}|\mathcal D|^{-\frac 12} e_p(-D\overline{n_1}s_1^2)\bigg|^2\bigg\}^{\frac 12}\\
\ge &\; \frac 12\; |\mathcal D|^{-\frac 12}\;\;\min_{D<s_1\le 2D}\bigg|\sum_{n_1\in\mathcal D}e_p(-D\overline{n_1}s_1^2)\bigg|.\end{aligned}
\ee

Note that by the definitions of $D$ and $\mathcal D$, $D\overline{n_1}\in\{1,\ldots, D\}$ for $n_1\in\mathcal D$ and $D\overline{n_1}s_1^2$ is an integer bounded by $p^{\frac 12}$ for $D<s_1<2D$. It follows that the inner sum in \eqref{2.28} is $\asymp |\mathcal D|$ and therefore \eqref{2.25} is
\be\label{2.29}
\gtrsim |\mathcal D|^{\frac 12} p^{-\frac 12}>\exp\bigg(c\;\frac{\log p}{\log\log p}\bigg)\; p^{-\frac 12}.\ee

Hence this example shows that the exponent $\frac 1{10}$ in \eqref{1.2} cannot be replaced by a clean $\frac 12$, though we are unable to rule out the validity of this inequality with exponent $\frac 12+\epsilon$.

\section
{Comments and another example}

The phenomenon described in Theorem 1.1 opens the door to a number of questions on possible generalizations
in different directions.

\noindent
{\bf Question 1.}  {\sl Let $\vp_1, \vp_2 \in\mathbb Z[X]$, $\vp_1(0)=\vp_2(0)=0$ be linearly independent
polynomials and define for $f_1, f_2:\mathbb F_p\to\mathbb C$
\be\label{3.1}
F(x) =\frac 1p \sum_{y\in\mathbb F_p} f_1\big(x+\vp_1(y)\big) f_2\big(x+\vp_2(y)\big).\ee
Is an inequality of the type \eqref{1.2} valid?}

\noindent {\bf Question 1'.}  {\sl Same question as above with $\vp_1, \vp_2$ a pair of rational functions, and excluding the poles in the summation \eqref{3.1}.
}

\medskip
Likely Theorem 1.1 and its proof extend to the case $\vp_1, \vp_2$ are linearly independent quadratic
polynomials.
Note that in the above analysis not only a bound
\be\label{3.2}
\Big|\sum_y e_p \big(a\vp_1(y)+ b\vp_2(y)\big)\Big| < c\sqrt p\ee
for $ab\not= 0$ (the latter being obtainable from Weil) is involved but also the exact evaluation of the
above exponential sum (which is possible only in special cases).

Beyond the bilinear case, one may also ask

\noindent{\bf Question 2.}
{\sl Is there a multilinear generalization of Theorem 1.1?
In particular if we define
\be\label{3.3}
F(x) =\frac 1p \sum_{y\in\mathbb F_p} f_1(x+y) f_2(x+y^2) f_3(x+y^3)\ee
does an inequality of the form
\be\label{3.4}
\big\Vert F-\mathbb E[f_1]\;\mathbb E[f_2] \;\mathbb E[f_3]\big\Vert_1\leq cp^{-\delta} \Vert f_1\Vert_\infty\,
\Vert f_2\Vert_\infty \,\Vert f_3\Vert_\infty.\ee
hold for some $\delta<0$?}
\medskip

The main theorem in \cite{FK} would imply that the left hand side of \eqref{3.4} is bounded by $\,o(1)\;\Vert f_1\Vert_\infty \; \Vert f_2\Vert_\infty \; \Vert f_3\Vert_\infty\,$
but without any quantitative specification (that would be awkward to extract from their approach) and certainly not a power gain $p^{-\delta}$ as in Theorem 1.1.
The explicit calculation based on the standard Fourier transform used above in the bilinear case does not seem to succeed for the trilinear average \eqref{3.3} that
likely would require higher order Fourier analysis.

Next, returning to Question 1', we establish another result in the spirit of Theorem 1.1, taking $\vp_1(x)=x, \vp_2(x)=\frac 1x$.
Thus we prove the following (with the same notations)

\noindent
{\bf Theorem 3.1.}
{\sl Define
\be\label{3.6}
F(x) =\frac 1p \sum_{y\in \mathbb F_p^*} f_1 (x+y) f_2\Big(x+\frac 1y\Big).\ee
Then the following inequality holds
\be\label{3.7}
\big\Vert F- \mathbb E [f_1]\; \mathbb E[f_2]\big\Vert_2 \leq cp^{-\frac 1{10}} \Vert f_1\Vert_2\cdot \Vert f_2\Vert_2.\ee}

Following the initial steps in the proof of Theorem 1.1, the left hand side of \eqref{3.7} is bounded by
\be\label{3.8}
\frac 1{\sqrt p} \bigg\{ \sum^{p-1}_{s=0} \bigg| \sum^{p-1}_{n=0} \hat f_{1}(n-s)\hat f_2 (n) K(s-n, n)\bigg|^2\bigg\}^{\frac 12}\ee
where now
\be\label{3.9}
K(x, y)=\begin{cases}
\frac 1{\sqrt p}\sum_{z\in \mathbb F_p^*} e_p (xz+y\frac 1z) \ & \text { if } \ y\not= 0\\
0 &\text { otherwise}
\end{cases}\ee
is given by a Kloosterman sum.

Note that by change of variables
\be\label{3.10}
\frac 1{\sqrt p} \sum_{z\in \mathbb F_p^*} e_p \Big(xz+y\frac 1z\Big) = \frac 1{\sqrt p} \sum_{z\in\mathbb F_p^*} e_p \Big(xyz+\frac 1z\Big)
= K\ell (xy)\ee
denoting the Kloosterman sum
\be\label{3.11}
\frac 1{\sqrt p} \sum_{z\in\mathbb F_p^*} e_p\Big(az+\frac 1z\Big) = K\ell_2 (a; p)=K\ell (a).\ee
(Note that $\overline{K\ell(a)}=K\ell(a)$.)

It is well-known that
\be\label{3.12}
|K\ell(a)|< c\ee
but $K\ell (a)$ does not allow an explicit evaluation.

Instead, we will rely on results from \cite {FKM} to bound $\Omega$ given by (2.13).
Thus, by (3.10)
\be\label{3.13}
\Omega =\Big\Vert \sum_x K\ell\big(x(s-x)\big) K\ell \big((x+u)(s-x)\big)\,K\ell\big(x(s'-x)\big) K\ell\big((x+u)(s'-x)\big)\Vert_{\ell^2_{s, s', u}}
\ee
and hence
\be\label{3.14}
\begin{aligned}\Omega^2 =&\sum_{x, y, s, s', u} K\ell\big(x(s-x)\big) \, K\ell\big((x+u)(s-x)\big) \, K\ell\big(x(s'-x)\big) \, K\ell \big((x+u)(s'-x)\big)\\
&\;\;\;\;\;\;\;\;\;\;\;\times K\ell \big(y(s-y)\big) \, {K\ell\big((y+u)(s-y)\big)} \,K\ell \big(y(s'-y)\big)\,  K\ell \big((y+u)(s'-y)\big).
\end{aligned} \ee

We will perform the summations in $s, s'$ with $x, y, u$ fixed.
Thus we bound
\be\label{3.15}
\sum^{p-1}_{s=0} K\ell \big(x(s-x)\big) K\ell\big((x+u)(s-x)\big) K\ell \big(y(s-y)\big) K\ell \big((y+u)(s-y)\big)\ee
which, in order to invoke \cite{FKM}, we rewrite as
\be\label{3.16}
\sum_sK\ell(\gamma_1.s)K\ell (\gamma_2.s) K\ell (\gamma_3. s) K\ell(\gamma_4.s)\ee
where $\gamma =\begin{pmatrix} a&b\\ c&d\end{pmatrix} \in PGL_2(\mathbb F_p)$ acts by linear fractional transformation
\be\label{3.17}
\gamma.s=\frac {as+b}{cs+d}\ee
and in our case \eqref{3.15} these transformations are affine
\be\label{3.18}
\begin{aligned}
\gamma_1 &= \begin{pmatrix} x&-x^2\\ 0&1\end{pmatrix} \qquad  \gamma_2 = \begin{pmatrix} x+u& -x(x+u)\\ 0&1\end{pmatrix}\\[14pt]
\gamma_3 &=\begin{pmatrix} y& -y^2\\ 0&1\end{pmatrix} \qquad \gamma_4= \begin{pmatrix} y+u& -y(y+u)\\ 0&1\end{pmatrix}
\end{aligned}
\ee
It follows in particular from Corollary 3.3 in \cite {FKM} that if $\gamma_1, \gamma_2, \gamma_3, \gamma_4\in PGL_2(\mathbb F_p)$ are pairwise distinct, then
\be\label{3.19}
\big|\sum_sK\ell(\gamma_1.s)K\ell (\gamma_2.s) K\ell (\gamma_3. s) K\ell(\gamma_4.s)\big|< c\sqrt p.\ee
By \eqref{3.19} the latter condition will be satisfied if $x\not=y$ and $u\not=0$.

Returning to \eqref{3.14}, we fix $x, y, u, x\not = y, u\not= 0$ and bound the $s, s'$ summations using \eqref{3.19}.
If $x=y$ or $u=0$, bound \eqref{3.16} trivially by $p$.
This gives a bound on \eqref{3.14}.
$$
\Omega^2< c (\sqrt p)^2 p^3+cp^2.p^2\leq cp^4.
$$
Therefore,
\be\label{3.20}
\Omega< cp^2.\ee
Substituting again the bound on $\Omega$ in \eqref{3.13} shows that
\be\label{3.21}
 \begin{aligned}&\frac 1{\sqrt p} \bigg\{ \sum^{p-1}_{s=0} \bigg| \sum^{p-1}_{n=0} \hat f_{1}(n-s)\hat f_2 (n) K(s-n, n)\bigg|^2\bigg\}^{\frac 12}\\\leq &\;\;c\frac 1{\sqrt p} p^{2/5} \Vert f_1\Vert_2 \cdot \Vert f_2\Vert_2 \leq cp^{-\frac 1{10}}\Vert f_1 \, \Vert_2 \dot \Vert f_2\Vert_2\end{aligned}
\ee
proving \eqref{3.7}.$\quad\square$

\end{document}